\documentstyle [12pt] {article}

\newenvironment{namelist}[1]{%
\begin{list}{}
    {
     
     \settowidth{\labelwidth}{#1}
     \setlength{\leftmargin}{1.1\labelwidth}
    }
  }{%
\end{list}}

\newcommand{\mapright}[1] {\stackrel{#1} {\hbox to 15pt{\rightarrowfill}}}
\newcommand{\maprightu}[1] {\stackrel{#1} {\hbox to 40pt{\rightarrowfill}}}
\newcommand{\mapdownr}[1] {\Big\downarrow
  \rlap{$\vcenter{\hbox{$\scriptstyle#1$}}$}}

\newcommand{\End}{\vskip-\baselineskip
  \hfill{\vrule width 6pt height 7.5pt depth 1.5pt}}

\newcommand{\Sym}[1]{\ifmmode{#1}\else \mbox{${#1}$}\fi}
\newcommand{\Vjp}{\bigbreak}
\newcommand{\ov}{\overline}

\newcommand{\wh}{\widehat}

\newcommand\setR{\mbox{\vrule height 7.2pt depth 0pt \hskip-0.8pt R}}

\newcommand\setZ{\bf Z}

\begin{document}

\baselineskip 18pt

\centerline{\large \bf A Natural Framing for Asymptotically Flat}

\centerline{\large \bf Integral Homology 3-Sphere}

\vskip 1cm
\centerline{\bf Su-Win Yang}

\vskip 10pt
\centerline{\it Communicated with W. H. Lin}

\vskip 1cm
\Vjp
\begin{abstract}
For an integral homology 3-sphere embedded asymptotically flatly
in an Euclidean space, we find a natural framing extending the standard
trivialization on the asymptotically flat part.
\end{abstract}

\vskip 12pt

Suppose $\ov M$ is a $3$-dimensional closed smooth manifold which has the same
integral homology groups as the $3$-sphere $S^3$. $x_0$ is a fixed point in
$\ov M$. Embed $\ov M$ in a Euclidean space $\setR^n$ such that $x_0$ is
the infinite point of the $3$-dimensional flat space $\setR^3 \times \{ 0 \}$
of $\setR^n$ and a neighborhood of $x_0$ contains the whole flat space
$\setR^3 \times \{ 0 \}$ except a compact set.
Precisely, for any positive number $r$, let $B_r$ denote the closed ball
of radius $r$ in $\setR^3$ and $N_r = ( \setR^3 - B_r ) \times \{ 0 \}$;
there exists $r_0$, a positive number, such that $N_{r_0}$ is contained
in $\ov M$ and $N_{r_0} \cup \{ x_0 \}$ is an open neighborhood of $x_0$
in $\ov M$.

Let $M = \ov M - \{ x_0 \}$, it is an asymptotically flat $3$-dimensional
manifold with acyclic homology. The main purpose of this article
is to define a natural framing for $M$.
If we identify the tangent spaces of points in the flat part $N_{r_0}$ with
$\setR^3 \times \{ 0 \}$, then the tangent bundle of $M$ can be thought as
a $3$-dimensional vector bundle over the closed manifold
$M_0 = M \big / \ov N_s$, where $s$ is a number greater than $r_0$ and $\ov N_s$
is the closure of $N_s$; we shall call this vector bundle the tangent bundle
$T(M_0)$ of $M_0$. And our natural framing is just a trivialization of $T(M_0)$,
which corresponds to a trivialization of the tangent bundle $T(M)$ whose restriction
to the flat part is the standard trivialization on $\setR^3$.
Because $M_0$ is a closed $3$-manifold, there are countably infinite many
choices of framings associated with the infinite elements in $[M_0, SO(3)]$.
( When $H_*(M_0) \approx H_*(S^3)$,
$[M_0, SO(3)] \approx [S^3, SO(3)] \approx \setZ$. )
Therefore, our natural framing is a special choice from the infinite many.

On the other hand, this natural framing for $T(M_0)$ can also provide a special
one-to-one correspondence between the infinite framings of $S^3$ and that of
$\ov M$. ( Note: Here, we do not think that $\ov M$ and $M_0$
have the same tangent
bundle. Conversely, we may think that the tangent bundle of $\ov M$ is equal to
the connected sum of the tangent bundles of $M_0$ and $S^3$. )

\vskip 20pt
There are two main steps to the natural framing on $T(M_0)$.

\vskip 10pt
\noindent
{\bf Step 1 \ \ A special map from $C_2(M)$ to $S^2$}

\vskip 8pt
We define $C_2(M)$ at first.

\vskip 6pt
For any set $X$, $\Delta(X)$ denote the diagonal subset
$\{ (x, x) \in X \times X, x \in X \}$ of $X\times X$ and
$C_2(X) = X \times X - \Delta(X)$.
Thus $C_2(M)$ is the configuration space of all pairs
of distinct two points in $M$.

\vskip 6pt
Fix some large number $s$ such that $M \subset (B_s \times \setR^{n-3})
\cup N_s$.

\vskip 6pt
For any $r \ge s$, let $B_r = \{ x \in \setR^3 : |x| \le r \}$,
$N_r = ( \setR^3 - B_r ) \times \{ 0 \}$ and $M_r = M - N_r$.

\vskip 6pt
Let $Y$ denote the union of the following three subsets of $C_2(M)$:
$$\mbox{(i)}\ \ Y_0 = C_2(N_s)\ \ \ \ $$
$$\mbox{(ii)}\ \ Y_1 = \cup_{r \ge s} (N_{r+s} \times M_r)$$
$$\mbox{(iii)}\ \ Y_2 = \cup_{r \ge s} (M_r \times N_{r+s})$$

\vskip 10pt
Let $\pi: \setR^n \longrightarrow \setR^3$ denote the projection
$$\pi(t_1, t_2, \cdots, t_n) = (t_1, t_2, t_3)$$
and $f: Y \longrightarrow S^2$ denote the map
$$f(x, y) = \frac{\pi(y - x)}{|\pi(y - x)|}$$
for $(x, y) \in Y$, $x, y \in M$.

For the well-defining of the map $f$, we should check that $|\pi(y - x)|$ is a
non-zero value. When $(x, y)$ is in $Y_0$, $|\pi(y - x)| = |y - x|$, it is
non-zero. When $(x, y)$ is in $Y_1$, $(x, y)$ is in $N_{r+s} \times M_r$
for some $r \ge s$; thus $\pi(x)$ is outside of $B_{r+s}$ and $\pi(y)$
is in $B_r$, and hence
$\pi(y - x) = \pi(y) - \pi(x)$, it has also a non-zero norm.
It is similar for the case that $(x, y)$ is in $Y_2$.

The following proposition describes some homology properties for
the space $Y$ and the map $f$.

\vskip 10pt
\noindent
{\bf Proposition 1}
\begin{namelist}{(iii)}
\item[{\rm (i)}]
$H_*(Y) \approx H_*(S^2)$
\item[{\rm (ii)}]
$f_*; H_2(Y) \longrightarrow H_2(S^2)$ is an isomorphism.
\item[{\rm (iii)}]
Let $j: Y \longrightarrow C_2(M)$ denote the inclusion map.
$$j_*: H_i(Y) \longrightarrow H_i(C_2(M))$$
is isomorphic, for all integer $i \ge 0$.
\End
\end{namelist}

In the proof of the proposition, we strongly use the assumption
that $H_*(M)$ is acyclic.

\vskip 8pt
\noindent
{\bf Remark:} \ All the homologies in this article are with
integral coefficients.

\vskip 9pt
By Proposition 1, the continuous map
$f: Y \longrightarrow S^2$ uniquely extends to
a continuous map $\ov f:C_2(M) \longrightarrow S^2$
up to homotopy relative to the subspace $Y$.
( That is, if both $\ov f_1$ and $\ov f_2$ are the extensions of $f$ to the
whole space $C_2(M)$, then there is a homotopy $F: C_2(M) \times [0,1]
\longrightarrow S^2$ such that $F(\xi, 0) = \ov f_1(\xi)$,
$F(\xi, 1) = \ov f_2(\xi)$, for all $\xi \in C_2(M)$, and $F(\xi', t) = f(\xi')$
for all $\xi' \in Y$ and $t \in [0, 1]$. )

Usually, the homotopy class of a map from $C_2(M)$ to $S^2$ can not give any
framing on $T(M_0)$. But the extension of $f$ does give a framing on $T(M_0)$
as shown in Step 2.

\vskip 20pt
\noindent
{\bf Step 2\ \ The framing determined by the map $\ov f$ on $C_2(M)$}

\vskip 10pt
The normal bundle of $\Delta(M)$ in $M \times M$ can be identified as
the tangent bundle $T(M)$ of $M$. Consider a suitable compactification of
$C_2(M)$, the spherical bundle $S(TM)$ become a part of boundary of $C_2(M)$.
Let $h: S(TM) \longrightarrow S^2$ denote the restriction of $\ov f$ to $S(TM)$.
On the flat part $N_s$ of $M$, the spherical bundle $S(TN_s) = N_s \times S^2$
and $h$ on $S(TN_s)$ is equal to the map restricted from $f$ which is exactly
the projection from $N_s \times S^2$ to $S^2$. Thus $h$ induces a map
$h_0: S(TM_0) \longrightarrow S^2$.

\vskip 6pt
$S(TM_0)$ is a $SO(3)$-bundle over $M_0$.

Can $h_0: S(TM_0) \longrightarrow S^2$ determine uniquely an orthogonal map,
that is, a fibrewise orthogonal map? ( An orthogonal map is exactly a framing
for the vector bundle. ) There is also an interesting question that can $h_0$
be homotopic to an orthogonal map; if such an orthogonal map exists,
is it unique up to homotopy?  We shall answer the questions partially.

Choose a framing for $S(TM_0)$ and we may think $h_0$ as a map from
$M_0 \times S^2$ to $S^2$. Let $y_0$ denote the point in $M_0$ representing
the set $N_s$. Then the restriction of $h_0$ to $y_0 \times S^2$ is the identity
map of $S^2$. Thus the restriction of $h_0$ to each fibre $x \times S^2$,
$x \in M_0$, is also a homotopy equivalence; and hence, $h_0$ induces a map
$\wh h_0$ from $M_0$ to $G(3)$, the space of all homotopy equivalences of $S^2$
to itself. Choose a base point $z_0$ in $S^2$, and consider the subspace $F(3)$
of $G(3)$ consisting of all the homotopy equivalences which fix the base point
$z_0$. Then $F(3)$ is the fibre of the fibration $G(3)$ over $S^2$, it is the
key fact for the homotopic computations.

\vskip 8pt
For any two spaces $X_1$ and $X_2$ with
base points $x_1$ and $x_2$, respectively,
$[X_1, X_2]$ denotes the set of
homotopy classes of continuous maps from $X_1$ to $X_2$
and sending $x_1$ to $x_2$. In the following, $M_0$is with base point $y_0$
representing the set $\ov N_s$; $SO(3)$, $G(3)$ and $F(3)$ are with the base
point the identity of $S^2$. We shall consider only the maps sending the base
point to base point and consider only the homotopies which keep the base point
fixed.

$M_0$ has the same homology as $S^3$. Usually, we can not expect they also
have the same homotopy behavior. But we still have the following proposition.

\vskip 10pt
\noindent
{\bf Proposition 2} \ \
Suppose $\phi: M_0 \longrightarrow S^3$ is a degree 1 map. Then the homotopy
classes $[M_0, SO(3)]$, $[M_0, G(3)]$, $[M_0, F(3)]$ are all groups, and
the group homomorphisms induced by $\phi$,
$$[S^3, SO(3)] \maprightu{\phi^{\sharp}} [M_0, SO(3)]$$
$$[S^3, G(3)] \maprightu{\phi^{\sharp}} [M_0, G(3)]$$
$$[S^3, F(3)] \maprightu{\phi^{\sharp}} [M_0, F(3)]$$
$$[S^3, S^2] \maprightu{\phi^{\sharp}} [M_0, S^2]$$
are all isomorphisms of groups.
\End

\vskip 12pt
There are further relations between these homotopy classes.

\vskip 8pt
\noindent
{\bf Proposition 3} \ \  Let $p: SO(3) \longrightarrow G(3)$ and
$q: F(3) \longrightarrow G(3)$ denote the inclusions.
Then, for any integral homology 3-sphere $M_0$, the homomorphism
$$p_* \oplus q_*: [M_0, SO(3)] \oplus [M_0, F(3)] \longrightarrow
[M_0, G(3)]$$
is an isomorphism.

\noindent
Especially, when $M_0 = S^3$, we have
$$\pi_3(G(3)) \approx \pi_3(SO(3)) \oplus \pi_3(F(3))\ .$$
\End

\vskip 10pt
Furthermore, the group isomorphism
$$q_*^{-1}: [M_0, G(3)] \big / p_*([M_0, SO(3)] \longrightarrow
[M_0, F(3)]$$
induces a group homomorphism
$$Q: [M_0, G(3)] \longrightarrow [M_0, F(3)] \approx \setZ_2 \ \ .$$

\vskip 6pt
For a continuous map $g: M_0 \times S^2 \longrightarrow S^2$,
let $\wh g$ denote the map from $M_0$ to $G(3)$ defined by
$\wh g(x)(y) = g(x, y)$, for $x \in M_0$ and $y \in S^2$
and let $Q(g) = Q([\wh g])$.

\vskip 8pt
\noindent
{\bf Theorem 4}\ \ A continuous map $g: M_0 \times S^2 \longrightarrow S^2$
is homotopic to an orthogonal map, if and only if, $Q(g) = 0$ in
$[M_0, F(3)]$.
\End

\vskip 8pt
Now, $h_0$ still denotes the map from $S(TM_0)$ to $S^2$ given by the map
$\ov f: C_2(M) \longrightarrow S^2$.
Choose a framing for $TM_0$, $\psi: S(TM_0) \longrightarrow M_0 \times S^2$,
it is a fibre map and fibrewise orthogonal.
Then $h_0 \circ \psi^{-1}$ is a map from $M_0 \times S^2$ to $S^2$
and the value $Q(h_0 \circ \psi^{-1})$ is independent of the choice of
the framing $\psi$. Therefore, $Q(h_0 \circ \psi^{-1})$ is an invariant
of the integral homology 3-sphere $\ov M$, it is the obstruction for $h_0$
to be homotopic to an orthogonal map. We hope that this is not really
an obstruction.

\vskip 10pt
\noindent
{\bf Conjecture 5} \ \ $Q(h_0 \circ \psi^{-1}) = 0$, for any integral
homology 3-sphere $\ov M$.
\End

\vskip 10pt
On the other hand, the group isomorphism
$$p_*^{-1}: [M_0, G(3)] \big / q_*([M_0, F(3)] \longrightarrow
[M_0, SO(3)]$$
induces a group homomorphism
$$P: [M_0, G(3)] \longrightarrow [M_0, SO(3)] \ \ .$$

\vskip 6pt
For a continuous map $g: M_0 \times S^2 \longrightarrow S^2$,
let $P(g) = P([\wh g])$.

\vskip 8pt
For the map $h_0$ and the corresponding element $P(h_0 \circ \psi^{-1})$ in
$[M_0, SO(3)]$, choose an orthogonal map
$g_0: M_0 \times S^2 \longrightarrow S^2$ such that the associated map $\wh g_0$
is in the homotopy class $P(h_0 \circ \psi^{-1})$.
Then we get an orthogonal map $g_0 \circ \psi: S(TM_0) \longrightarrow S^2$
which represents a homotopy class of framings determined by $h_0$, also by the
map $\ov f: C_2(M) \longrightarrow S^2$. This framing can also be characterized
by the following theorem.

\vskip 10pt
\noindent
{\bf Theorem 6}\ \ There exists a framing $\psi_0: S(TM_0) \longrightarrow
M_0 \times S^2$ unique up to homotopy such that
$P(h_0 \circ \psi_0^{-1}) = 0$.
\End

\vskip 20pt
\noindent
{\bf Proofs}

\vskip 10pt
\noindent
{\bf Outline of Proof of Proposition 1}

$N_s$ is a subset of $\setR^3 \times \{ 0 \}$.
In $N_s$, we choose a subspace $S_3$ which is a deformation retract of $N_s$
and a point $x_1$ in the bounded component of $\setR^3 \times \{ 0 \} - S_3$.
Let $S = \{ x_1\} \times S_3$, it is a subspace of $Y$.
We show that the three maps, the inclusion of $S$ in $Y$, the restriction of
$f$ to $S$, and the restriction of $j$ to $S$, all induce isomorphisms of
homology groups of the corresponding spaces.
That is, $H_*(S) \longrightarrow H_*(Y)$,
$(f|_S)_*: H_*(S) \longrightarrow H_*(S^2)$,
and $(j|_S)_*: H_*(S) \longrightarrow H_*(C_2(M))$ all are isomorphisms.

\vskip 10pt
\noindent
{\bf Proof of Proposition 1}

\vskip 8pt
First we compute the homology of $Y_0, Y_1, Y_2$, separately.

\vskip 6pt
$Y_0 = C_2(N_s) = N_s \times N_s - \Delta(N_s) \subset N_s \times N_s$.
$N_s$ is homeomorphic to $S^2 \times (s, \infty)$.
Thus $H_*(N_s \times N_s) \approx H_*(S^2 \times S^2)$.
By Thom Isomorphism, $H_i(N_s \times N_s, Y_0) \approx H_{i-3}(N_s)$.

Now, we use the long exact sequence of the pair $(N_s \times N_s, Y_0)$ to
determine $H_*(Y_0)$.

\vskip 8pt
$\longrightarrow H_{i+1}(N_s \times N_s, Y_0) \maprightu{\partial_*} H_i(Y_0)
\longrightarrow H_i(N_s \times N_s) \longrightarrow$

$\longrightarrow$ $H_i(N_s \times N_s, Y_0) \longrightarrow \cdots$

\vskip 8pt
When $i$ is odd, both $H_{i+1}(N_s \times N_s, Y_0)$ and $H_i(N_s \times N_s)$
are the trivial group $\{ 0 \}$. Thus we have
$$H_4(Y_0) \approx H_4(N_s \times N_s) \oplus
\partial_*(H_5(N_s \times N_s, Y_0) \approx \setZ \oplus \setZ$$
$$H_2(Y_0) \approx H_2(N_s \times N_s) \oplus
\partial_*(H_3(N_s \times N_s, Y_0)
\approx \setZ \oplus \setZ \oplus \setZ$$
and $H_i(Y_0)$ is trivial, if $i$ is odd.

( $\setZ$ denotes the group of integers. )

\vskip 8pt
To find the generators of $H_4$ and $H_2$ of $Y_0$, we choose three 2-spheres
$S_1, S_2, S_3$ in $\setR^3 \times \{ 0 \}$ of radius $2s, 4s, 6s$,
respectively, all with center the origin.
( $S_i$ is the boundary of $N_{2s \times i}$, $i = 1, 2, 3$. )
For each $i$, $i = 1, 2, 3,$ choose a point $x_i$ in $S_i$. The 2-spheres are
also oriented in the same way, that is, the natural diffeomorphisms of the
2-spheres are orientation-preserving.
Then $S_i \times x_j$ and $x_j \times S_i$, $1 \le i \ne j \le 3$, are 2-cycles
in $Y_0$, also in $N_s \times N_s$;
$S_i \times S_j$, $1 \le i \ne j \le 3$, are 4-cycles
in $Y_0$, also in $N_s \times N_s$.

In the following, if $c$ is a cycle in $Y_0$,
$[c]$ shall denote the corresponding
homology class in $Y_0$.

\vskip 8pt
\noindent
{\bf Lemma 7}

(i) $[S_1 \times S_3]$ is the generator of $H_4(N_s \times N_s)$.

(ii) $[(S_1 - S_3) \times S_2]$ is the generator of the subgroup
$\partial_*(H_5(N_s \times N_s, Y_0))$ in $H_4(Y_0)$.
\End

\vskip 8pt
We use the lemma to prove Proposition 1, and prove the lemma later.

\vskip 8pt
There are some relations between these classes in $H_*(Y_0)$:\\
$[(S_1 - S_3) \times S_2] = [S_1 \times S_2] - [S_3 \times S_2]$,\ \ \
$[(S_1 \times S_2] = [S_1 \times S_3]$\\
and $[S_3 \times S_2] = [S_3 \times S_1]$.

\vskip 6pt
Thus $[S_1 \times S_3]$ and $[S_3 \times S_1]$
form the basis of $H_4(Y_0)$.

Similarly, $[S_1 \times x_3]$ and $[x_1 \times S_3]$ are the basis of
$H_2(N_s \times N_s)$;\\
$[(S_1 - S_3) \times x_2]$
( $= \epsilon_0 [x_2 \times (S_1 - S_3)]$, $\epsilon_0$ is 1 or $-1$ )
is the generator of the subgroup
$\partial_*(H_3(N_s \times N_s, Y_0)$ in $H_2(Y_0)$.

\vskip 6pt
Thus $[S_1 \times x_3]$, $[x_1 \times S_3]$
and $[(S_1 - S_3) \times x_2]$ form a basis of $H_2(Y_0)$.

\vskip 8pt
Now we study the homology of $Y_1$ and $Y_2$.

It is easy to see that the inclusion of $N_{4s} \times M_{3s}$ in $Y_1$
and the inclusion of $S_3 \times M_{3s}$ in $N_{4s} \times M_{3s}$
both are homotopy equivalences. Thus $H_*(Y_1) \approx H_*(S_3 \times M_{3s})$
$\approx H_*(S_3)$. ( Recall: $M_r$ is acyclic, for any $r \ge s$. )
Similarly, $Y_2$ also has the same homology as 2-sphere.

$Y_1$ and $Y_2$ are disjoint, and hence the homology
of their union $Y_1 \cup Y_2$ is also determined.
We can use the Mayer-Vietoris Sequence of the
triple $(Y, Y_0, Y_1 \cup Y_2)$ to find the homology of $Y$. In fact, we have
\begin{namelist}{(iv)}
\item[{\rm (i)}]
The cycle $S_1 \times S_3$ is contained in $Y_2$ and is killed in $Y_2$.
\item[{\rm (ii)}]
The cycle $S_3 \times S_1$ is contained in $Y_1$ and is killed in $Y_1$.
\item[{\rm (iii)}]
The cycle $x_3 \times S_1$ is contained in $Y_1$ and is killed in $Y_1$.
\item[{\rm (iv)}]
The cycle $S_1 \times x_3$ is contained in $Y_2$ and is killed in $Y_2$.
\end{namelist}

Therefore, $H_4(Y) = \{ 0 \}$ and in $H_2(Y)$,
we have $[x_1 \times S_3]$
and $[S_3 \times x_2]$ left;
the equality
$[(S_1 - S_3) \times x_2]
= \epsilon_0 [x_2 \times (S_1 - S_3)]$
become the new equality $-[S_3 \times x_2] = - \epsilon_0 [x_2 \times S_3]$.
Thus $[x_1 \times S_3] = [x_2 \times S_3]$
$= \epsilon_0 [S_3 \times x_2] = \epsilon_0 [S_3 \times x_1]$.
This proves that $H_*(Y) \approx H_*(S^2)$. Actually, we know more than that:
the inclusion of the space $\{ x_1 \} \times S_3$ in $Y$ induces isomorphisms
of the homology groups. It is easy to see that the map $f$, restricted to
$\{ x_1 \} \times S_3$, is an homotopy equivalence from $\{ x_1 \} \times S_3$
to $S^2$. This proves the second statement that $f_*$ is an isomorphism.

To prove the third statement that $j_*$ is an isomorphism, it is also enough
to show that the restriction of $j$ to $\{ x_1 \} \times S_3$ induces
isomorphisms for the homology groups.
Similar to the computation of the homology of $C_2(N_s)$,
we consider the long exact sequence of pair $(M \times M, C_2(M))$

\vskip 6pt
$H_{i+1}(M \times M)
\longrightarrow H_{i+1}(M \times M, C_2(M)) \maprightu{\partial_*} H_i(C_2(M))$

$\longrightarrow
H_i(M \times M)$ $\longrightarrow$.

\vskip 8pt
Because $H_*(M)$ is acyclic, $H_*(M \times M)$ is also acyclic.

We have
$$H_i(C_2(M)) \approx
H_{i+1}(M \times M, C_2(M)), \  \mbox{for all} \  i \ge 1 \ .$$
But $H_{i+1}(M \times M, C_2(M)) \approx H_{i-2}(M)$, by the Thom Isomorphism.
Thus $C_2(M)$ has the same homology as 2-sphere.

And it is easy to see that the inclusion of $\{ x_1 \} \times (M, M - x_1)$
in $(M \times M, C_2(M))$ induces isomorphisms of homology groups, and hence,
the inclusion of $\{ x_1 \} \times (M - x_1)$ in $C_2(M)$ also induces
isomorphisms of homology groups.
The cycle $\{ x_1 \} \times S_3$ is a generator of
$H_2(\{ x_1 \} \times (M - x_1))$,
and hence also a generator of $H_2(C_2(M))$.
This proves the third statement that $j_*$ is an isomorphism.

\vskip 8pt
(i) of Lemma 7 is obvious. Now, we are going to prove (ii) in Lemma 7.
Consider the following commutative diagram
$$\matrix{
H_2(S_2) \otimes H_3(N_s, N_s - S_2) & \maprightu{id \otimes \partial_*}
& H_2(S_2) \otimes H_2(N_s - S_2) \cr\cr
\mapdownr{\tau_1} & & \mapdownr{\eta_1} \cr\cr
H_5(S_2 \times N_s, S_2 \times N_s - S_2 \times S_2)  & \maprightu{\partial_*}
& H_4(S_2 \times N_s - S_2 \times S_2) \cr\cr
\mapdownr{\tau_2} & & \mapdownr{\eta_2} \cr\cr
H_5(S_2 \times N_s, S_2 \times N_s - \Delta(S_2))  & \maprightu{\partial_*}
& H_4(S_2 \times N_s - \Delta(S_2)) \cr\cr
\mapdownr{\tau_3} & & \mapdownr{\eta_3} \cr\cr
H_5(N_s \times N_s, Y_0) & \maprightu{\partial_*} & H_4(Y_0) \cr\cr}
$$

The maps $\tau_1$ and $\eta_1$ are isomorphisms from Kunneth formula.
Other homomorphisms are induced by the corresponding inclusion maps.
$\tau_2$ is an isomorphism by the result of Lefschetz Duality in the
$5$-dimensional manifold $S_2 \times N_s$; $\tau_3$ is an isomorphism
by the result of Thom Isomorphism Theorem. Precisely, consider the following
commutative diagram
$$\matrix{
H_5(S_2 \times N_s, S_2 \times N_s - S_2 \times S_2)  & \maprightu{\sigma_1}
& H^0(S_2 \times S_2) \cr\cr
\mapdownr{\tau_2} & & \mapdownr{\tau_4} \cr\cr
H_5(S_2 \times N_s, S_2 \times N_s - \Delta(S_2))  & \maprightu{\sigma_2}
& H^0(\Delta(S_2)) \cr\cr}
$$
where $\sigma_i, i =1, 2,$ are the isomorphisms of Lefschetz Duality,
$\tau_4$ is the homomorphism induced by the inclusion.

\vskip 6pt
Because $\tau_4$ is an isomorphism, $\tau_2$ is also an isomorphism.
The proof of isomorphism of $\tau_3$ is in some sense
analogous to that for $\tau_2$, we omit it.

>From the long exact sequence of the pair $(N_s, N_s - S_2)$,
it is easy to see that $[S_1 - S_3]$ is the generator of
$\partial_*(H_3(N_s, N_s - S_2))$, and hence, $[S_2 \times (S_1 - S_3)]$
is the generator of
$(id \otimes \partial_*)(H_2(S_2) \otimes H_3(N_s, N_s - S_2))$.
By the commutativity of the above diagram, $[S_2 \times (S_1 - S_3)]$
($= - [(S_1 - S_3) \times S_2]$)
is the generator of $\partial_*(H_5(N_s \times N_s, Y_0))$.
This proves Lemma 7
and completes the long proof of {\bf Proposition 1}.

\vskip 10pt
\noindent
{\bf Proof of Proposition 2}

We need to show the isomorphisms between $[M_0, X]$ and $[S^3, X]$,
for $X = SO(3), G(3), F(3)$ and $S^2$.

For the case of $SO(3)$, we consider the classifying space $BSO(3)$ of the
$SO(3)$-bundles. Then
$$[M_0, SO(3)] \approx [SM_0, BSO(3)]\ \ \mbox{and}\ \
[S^3, SO(3)] \approx [S^4, BSO(3)]\ ,$$
where $SM_0$ is the suspension of $M_0$.
On the other hand, because $SM_0$ is simply connected and the map $S(\phi):
SM_0 \longrightarrow SS^3 (= S^4)$ induces isomorphisms of homology groups,
$S(\phi)$ is a homotopy equivalence. Thus $S(\phi)^{\sharp}:$$[SM_0, BSO(3)]$
$\longrightarrow [S^4, BSO(3)]$ is isomorphic, and hence,
$$[M_0, SO(3)] \approx [S^3, SO(3)]\ \ .$$

For the cases of $G(3)$ and $F(3)$, we may also consider the corresponding
classifying spaces, by the result of Fuchs [2];
and the proof is completely similar.

The group property of the associated homotopy classes is
a result of Dold and Lashof [1]; for the convenience of interested reader,
we give a proof in the appendix.

For the case of $S^2$, it is enough to note that
$[M_0, S^2] \approx [M_0, S^3]$ ( $\approx H^3(M_0)$ ), which implies the
isomorphism we need.

\vskip 10pt
\noindent
{\bf Proof of Proposition 3}

By Proposition 2, it is enough to prove the result for the case that
$M_0 = S^3$.

Consider the commutative diagram of fibrations over $S^2$
$$\matrix{
S^1 & \maprightu{} & SO(3) & \maprightu{} & S^2 &\cr\cr
\mapdownr{} & & \mapdownr{p} & & \mapdownr{id} & \cr\cr
F(3) & \maprightu{q} & G(3) & \maprightu{} & S^2 & \cr\cr}
$$
and the associated commutative diagram of exact sequences of homotopy groups
$$\matrix{
\pi_i(S^1) & \maprightu{} & \pi_i(SO(3)) & \maprightu{} & \pi_i(S^2) &\cr\cr
\mapdownr{} & & \mapdownr{p_*} & & \mapdownr{id} & \cr\cr
\pi_i(F(3)) & \maprightu{q_*} & \pi_i(G(3)) & \maprightu{\alpha} & \pi_i(S^2) & \cr\cr}
$$

For $i \ge 3$, $\pi_i(S^1) = \pi_{i-1}(S^1) = \{ 0 \}$, and hence
$$\pi_i(SO(3)) \approx \pi_i(S^2) \ \ .$$
Thus $p_*: \pi_i(SO(3)) \longrightarrow \pi_i(G(3))$ can be thought as the
right-inverse of $\alpha: \pi_i(G(3)) \longrightarrow \pi_i(S^2)$.
This implies that $\alpha$ is an epimorphism, $q_*$ is a monomorphism, and
$p_*$ supplies the necessary homomorphism for splitting.
Therefore, $$\pi_i(G(3)) = q_*(\pi_i(F(3)) \oplus p_*(\pi_i(SO(3)),
\  \mbox{for all}\ i\ge3$$

\newpage
\noindent
{\bf Appendix}

\vskip 8pt
The proof of the appendix is essentially from the proof of the main result
in Dold and Lashof [1]. The author just write it for self-interesting.

\vskip 10pt
Suppose $H$ is a path-connected space and has an associative multiplication
which  has a two-sided unit $e$. For $h_1, h_2 \in H$,
$h_1 \cdot h_2$ denotes the
product of $h_1$ and $h_2$. Thus $h \cdot e = e \cdot h = h$, for all $h \in H$.
Furthermore, assume $X$ is a polyhedron. The purpose of this appendix is
to show that the homotopy classes in $[X, H]$ form a group under the following
multiplication:

For any two maps $f, g: X \longrightarrow H$, $(f \cdot g)(x) = f(x) \cdot g(x)$.

\vskip 6pt
The associative law of this multiplication in $[X, H]$ is obvious. It is enough
to show that for any $f: X \longrightarrow H$, there is a map
$g: X \longrightarrow H$ such that $f \cdot g$ is homotopic to the constant map
$\ov e: X \longrightarrow H$, $\ov e(x) = e$, for all $x \in X$.

\vskip 8pt
We shall construct the map $g: X \longrightarrow H$ and the homotopy
$D: X \times I \longrightarrow H$ satisfying $D(x, 0) = e$,
$D(x, 1) = f(x) \cdot g(x)$, inductively on the skeleton of $X$.
( $I$ is the unit interval $[0, 1]$. )

$X^{(k)}$ denotes the $k$-skeleton of $X$.

Assume $g$ is defined on $X^{(k)}$ and $D$ is defined on $X^{(k)} \times I$
such that $D(x, 0) = e$ and $D(x, 1) = f(x) \cdot g(x)$,
for all $x \in X^{(k)}$. If necessary, we may ask that the base point
$x_0$ of $X$ is in $X^{(0)}$ and $f(x_0) = g(x_0)$ $=$ $D(x_0, t) = e$,
for any $t \in I$.

For any $(k+1)$-simplex $\Delta$ in $X^{(k+1)}$, we want to extend $g$
to the part $\Delta$ and $D$ to the part $\Delta \times I$. Let $S$ denote
the boundary of $\Delta$, it is a $k$-sphere. $S$ is in $X^{(k)}$, $g$ is
defined on $S$ and $D$ is defined on $S \times I$.

\vskip 8pt
\noindent
{\bf Claim}\ \ $g|_S: S \longrightarrow H$ is null-homotopic.

{\bf Proof}\ \ $\Delta$ is a simplex, there is a contraction map $\gamma:
\Delta \times I \longrightarrow \Delta$, $\gamma(x, 0) = x$ and $\gamma(x, 1) =
x_1$, for all $x \in \Delta$. $x_1$ is some fixed point in $S$.
Let $\beta: S \times I \longrightarrow H$ denote the map
$\beta(x, t) = f(\gamma(x, t)) \cdot g(x),$ for $x \in S$. Let $y_1 = f(x_1)$
and $\ov y_1: S \longrightarrow H$ denote the constant map sending the points
of $S$ to $y_1$. Then $\beta$ is a homotopy between $f \cdot g$ and $\ov y_1
\cdot g$ on $S$. $H$ is path-connected, $\ov y_1 \cdot g$ is homotopic to
$\ov e \cdot g = g$. Thus $g$ is homotopic to $f \cdot g$ on $S$.
On the other hand,
the restriction of $D$ to $S \times I$ provides a homotopy between
the restrictions of $f \cdot g$ and $\ov e$.
This proves that $g|_S$ is null-homotopic.

\vskip 8pt
Therefore, we can extend $g|_S$ to the part $\Delta$,
say, $g': \Delta \longrightarrow H$, and we can also extend $D|_{S \times I}$
to the whole boundary of $\Delta \times I$ as follows:
\subitem
We use $D': \partial (\Delta \times I) \longrightarrow H$ to denote
the extension. $\partial (\Delta \times I) = \Delta \times \{ 0 \}$ $\cup
\Delta \times \{ 1 \}$$\cup S \times I$.\\
$D'(x, 0) = e$ and $D'(x, 1) = f(x) \cdot g'(x)$, for all $x \in \Delta$; \\
$D'(y, t) = D(y, t)$, for all $y \in S$ and $t \in I$.

The map $D'$ may not be extended to $\Delta \times I$.
We shall find a map $g_1: \Delta \longrightarrow H$ with $g_1|_S = \ov e|_S$
and modify the map $D'$ by multiplying $D'$ with $g_1$ on the
part $\Delta \times \{ 1 \}$ such that the new map is null-homotopic.
Precisely, let $D'': \partial (\Delta \times I) \longrightarrow H$
denote the map, $D''(\xi) = D'(\xi)$, for all $\xi \in
\Delta \times \{ 0 \}$ $\cup S \times I$, $D''(x, 1) = D'(x, 1) \cdot g_1(x)$,
for all $(x, 1) \in \Delta \times \{ 1 \}$.

We may think the map $g_1$ as a map on $\Delta \times \{ 1 \}$
and extend it trivially to the whole boundary
$\partial (\Delta \times I)$, that is, sending all points undefined to $e$.
Then $D''$ is just equal to $D' \cdot g_1$.
To let $D''$ be null-homotopic, we can choose $g_1$
such that $[g_1] = [D']^{-1}$ in $\pi_{k+1}(H)$.
Of course, $g'$ should be changed to the new map $g' \cdot g_1$.
Therefore, $D''$ is null-homotopic and its extension to $\Delta \times I$
also gives the homotopy between $f \cdot (g' \cdot g_1)$ on $\Delta$.
This finishes the extension of $g$ and $D$ to $\Delta$.

\vskip 20pt
\centerline{\Large \bf References}

\begin{namelist}{[1]}
\item[{\rm [1]}]
A. Dold and R. Lashof, {\it Principal quasi-fibrations and fibre
homotopy equivalence of bundles}, Illi. J. Math. {\bf 3} (1959), 285-305.
\item[{\rm [2]}]
M. Fuchs, {\it, Verallgemeinerte Homotopie-Homomorphismen und
klassifizierende Raume}, Math. Ann. {\bf 161}, (1965), 197-230.
\end{namelist}

\noindent{Department of Mathematics} \\
\noindent{National Taiwan University} \\
\noindent{Taipei, Taiwan} \\
\noindent{E-mail: swyang@math.ntu.edu.tw}

\end{document}